\documentclass[12pt]{article}
\usepackage{amssymb,latexsym}
\usepackage{graphicx}
\def\ZZ{{\mathbb Z}}

\def\RR{{\mathbb R}}

\newtheorem{formula}{}[section]

\newtheorem{definition}[formula]{\indent Definition}
\newtheorem{corollary}[formula]{\indent Corollary}
\newtheorem{remark}[formula]{\indent Remark}
\newtheorem{lemma}[formula]{\indent Lemma}
\newtheorem{theorem}[formula]{\indent Theorem}

\def\thrm{\begin{theorem}}
\def\thrml#1{\begin{theorem}\label{#1}}
\def\ethrm{\end{theorem}}
\def\rmrk{\begin{remark}}
\def\rmrkl#1{\begin{remark}\label{#1}}
\def\ermrk{\end{remark}}
\def\dfntn{\begin{definition}}
\def\dfntnl#1{\begin{definition}\label{#1}}
\def\edfntn{\end{definition}}
\def\nmrt{\begin{enumerate}}
\def\enmrt{\end{enumerate}}

\def\qtn{\begin{equation}}
\def\qtnl#1{\begin{equation}\label{#1}}
\def\eqtn{\end{equation}}
\def\lmm{\begin{lemma}}
\def\lmml#1{\begin{lemma}\label{#1}}
\def\elmm{\end{lemma}}
\def\crllr{\begin{corollary}}
\def\crllrl#1{\begin{corollary}\label{#1}}
\def\ecrllr{\end{corollary}}

\begin{document}
\title{}
\date{}
\maketitle
\vspace{-0,1cm} \centerline{\bf Entropy of radical ideal of a tropical prevariety}
\vspace{7mm}
\author{
\centerline{Dima Grigoriev}
\vspace{3mm}
\centerline{CNRS, Math\'ematique, Universit\'e de Lille, Villeneuve
d'Ascq, 59655, France} \vspace{1mm} \centerline{e-mail:\
dmitry.grigoryev@univ-lille.fr } \vspace{1mm}
\centerline{URL:\ http://en.wikipedia.org/wiki/Dima\_Grigoriev} }

\begin{abstract}
The entropy of a tropical ideal is introduced. The radical of a tropical ideal consists of all tropical polynomials vanishing on
the tropical prevariety determined by the ideal. We prove that the entropy of the radical of a tropical bivariate polynomial with
vanishing coefficients equals zero.  Also we prove that the entropy of a zero-dimensional tropical prevariety vanishes. An example of a non-radical tropical ideal having a positive entropy is exhibited.
\end{abstract}

{\bf keywords}: radical of a tropical ideal, entropy


\section*{Introduction}

One can find the basic concepts of tropical mathematics in \cite{MS}.

Let $f=\min_{1\le j\le m} \{a_j+\sum_{1\le i\le n} t_{j,i}X_i\}$ be a tropical polynomial. Consider a family of {\it linearizations}
of $f$:
$$\min_{1\le j\le m} \{a_j + u(t_{j,1}+s_1,\dots,t_{j,n}+s_n)\}$$
\noindent in $N^n$  variables $u(k_1,\dots,k_n),\, 0\le k_1,\dots,k_n < N$ for some $N$ where $s_1,\dots,s_n \in \ZZ$,
provided that $0\le t_{j,1}+s_1,\dots,t_{j,n}+s_n < N$. Observe that if a point $(x_1,\dots,x_n)\in \RR^n$ satisfies $f$ (i.~e. the
minimum in $f$ is attained at least twice \cite{MS}) then the point 
$$u:=\{u(k_1,\dots,k_n)=k_1x_1+\cdots +k_nx_n\, : \,
0\le k_1,\dots,k_n < N\} \in \RR^{N^n}$$
\noindent satisfies the linearizations of $f$. Denote by $U_N\subset \RR^{N^n}$ a tropical linear prevariety \cite{MS} of 
the points  satisfying all the linearizations of $f$. 

We establish existence of the limit
$$H:=H(f):=\lim_{N\to \infty} \dim(U_N)/N^n$$
\noindent and call it the {\it (tropical) entropy} of $f$. Evidently, $0\le H\le 1$. In the univariate case ($n=1$) the tropical
entropy was introduced and studied in \cite{G}. In a similar way one extends the definition of the entropy $H(I)$ to tropical ideals $I$.

Informally speaking, the tropical entropy plays a role similar to the coefficient at $n$-th power of Hilbert's polynomial of an
ideal. In the classical commutative algebra this coefficient obviously vanishes for any non-zero ideal (Hilbert's polynomial has
the degree at most $n-1$ being equal the dimension of the variety determined by the ideal). This is not the case in the tropical
setting: the tropical entropy can be positive, an example of this phenomenon is provided in section~\ref{three}.

For a tropical ideal $I$ its tropical prevariety $V(I)\subset \RR^n$ consists of all tropical solutions of $I$ (recall that $V(I)$ is a
finite union of convex polyhedra \cite{MS}). We define in section~\ref{one} the {\it radical} $rad(I)$ of a tropical ideal $I$ as the set of all tropical polynomials
vanishing on $V(I)$. Unlike Hilbert's strong Nullstellensatz which describes the radical of an ideal in the classical commutative
algebra, the structure of the radical of a tropical ideal is more complicated. We mention also that a tropical version of Hilbert's
weak Nullstellensatz was obtained in \cite{GP}.

The main result of section~\ref{two} states that for a tropical bivariate polynomial $f:=\min_{1\le j\le m} \{t_{j,1}X+t_{j,2}Y\}$ with zero coefficients (so, whose prevariety is a tropical curve with a single vertex)  the entropy of its radical
$H(rad(f))=0$ vanishes. In \cite{G} it was proved in the univariate case ($n=1$) that $H(f)=0$ iff $f=rad(f)$ and moreover,
when $H(f)>0$ it holds $H(f)\ge 1/6$.

We prove in section~\ref{four} that the entropy of the radical of a zero-dimensional tropical prevariety (so, of a finite number of points) equals zero.

It would be interesting to clarify, whether one can generalize both results of this paper to vanishing the entropy of the radical of an arbitrary tropical ideal.

\section{Radical of a tropical ideal and entropy}\label{one}

Consider a tropical polynomial
\begin{eqnarray}\label{1}
f=\min_{1\le j\le m} \{a_j+\sum_{1\le i\le n} t_{j,i}X_i\}
\end{eqnarray}
where $a_j+\sum_{1\le i\le n} t_{j,i}X_i$ being linear functions (tropical monomials) with integers $t_{j,i},\, 1\le i\le n$ and
$a_j\in \RR$.  A point $(x_1,\dots,x_n)\in \RR^n$ is a tropical solution of $f$ if the minimum in (\ref{1}) is attained at least
twice \cite{MS}. The set of all tropical solutions of $f$ is called the tropical prevariety $V(f)\subset \RR^n$ of $f$. More generally,
one defines the tropical prevariety $V(I)\subset \RR^n$ of a tropical ideal $I$.

We define the radical $rad(I)$ of $I$ as the set (a tropical ideal) of all tropical polynomials vanishing on $V(I)$. Unlike Hilbert's
strong Nullstellensatz the radical of a tropical ideal is not exhausted by extracting roots of elements of the ideal (we'll see some
examples below). We mention that a tropical version (in a dual form) of Hilbert's weak Nullstellensatz was established in \cite{GP}.
\vspace{1mm}

In \cite{G} the entropy of a tropical polynomial $f$ was introduced as follows. For an integer $N$ consider a tropical prevariety
$U_N\subset \RR^{N^n}$ consisting of points $\{u(k_1,\dots,k_n)\in \RR\, :\, 0\le k_1,\dots,k_n< N\}$ satisfying tropical linear equations
\begin{eqnarray}\label{3}
\min_{1\le j\le m} \{a_j+u(t_{j,1}+s_1,\dots,t_{j,n}+s_n)\}
\end{eqnarray}  
over the variables $\{u(k_1,\dots,k_n)\, :\, 0\le k_1,\dots,k_n< N\}$ for any vector $(s_1,\dots,s_n)\in \ZZ^n$, provided that
$0\le t_{j,1}+s_0,\dots,t_{j,n}+s_n < N,\, 1\le j\le m$. We call (\ref{3}) the {\it linearization} of the tropical polynomial
\begin{eqnarray}\label{9}
\min_{1\le j\le m} \{a_j+\sum_{1\le i\le n} (t_{j,i}+s_i)X_i\}.
\end{eqnarray}

Observe that if a point $x:=(x_1,\dots,x_n)\in \RR^n$ is a solution of $f$ then $x$ satisfies also (\ref{9}) and the point
$$\{u(k_1,\dots,k_n)=k_1x_1+\cdots+k_nx_n\, :\, 0\le k_1,\dots,k_n< N\} \in U_N$$
\noindent due to (\ref{1}), (\ref{3}), (\ref{9}). On the other hand, $U_N$ can contain points not arising from tropical solution of $f$ (we'll see examples below).

Consider a partition of $n$-dimensional  grid $T_N:=\{(k_1,\dots,k_n)\, :\, 0\le k_1,\dots,k_n< N\} \subset \ZZ^n$
with the side $N$ into subgrids with sides $q_1,\dots,q_R$, respectively. Then the number of points in $T_N$ equals
$N^n=q_1^n+\cdots+q_R^n$. Denote by $p_r\,:\, \RR^{N^n} \twoheadrightarrow \RR^{q_r^n},\, 1\le r\le R$ the
projection of the coordinates from $T_N$ onto the coordinates from the $r$-th subgrid. Then $p_r(U_N)\subset U_{q_r}$ and
$U_N\subset U_{q_1}\times \cdots \times U_{q_R}$. Hence 
\begin{eqnarray}\label{4}
\dim(U_N)\le \dim(U_{q_1})+\cdots +\dim(U_{q_R}).
\end{eqnarray}

Therefore, similar to the proof of Fekete's subaddivity lemma \cite{S} one can verify that there exists a limit
\begin{eqnarray}\label{5}
H:=H(f):=\lim_{N\to \infty} \dim (U_N)/N^n =\inf \dim (U_N)/N^n.
\end{eqnarray}
Indeed, for any fixed $q$ partition grid $T_N$ into $\lfloor N/q\rfloor ^n$ subgrids equal $T_q$ which fill grid $T_{q\cdot \lfloor N/q\rfloor} \subset T_N$ and $N^n - (q\cdot \lfloor N/q\rfloor)^n = O(N^{n-1})$ subgrids each equal $T_1$. Due to (\ref{4})
$$\dim (U_N)\le \lfloor N/q \rfloor ^n \cdot \dim (U_q) +  N^n - (q\cdot \lfloor N/q\rfloor)^n.$$
\noindent With $N$ tending to the infinity, we conclude that 
$$\lim \sup_{N\to \infty}  \dim (U_N)/N^n \le \dim (U_q)/q^n$$
\noindent which implies (\ref{5}).

Note that in \cite{G} the entropy was defined by means of considering parallelepipeds (rather than cubes as in the present paper).
One can verify that these two definitions of the entropy coincide.

We call $H(f)$ the {\it (tropical) entropy} of $f$. More generally, one defines in a similar way the entropy $H(I)$ of a tropical ideal $I$. Clearly, if $I\subset I_1$ then $H(I)\ge H(I_1)$. Obsiously, $0\le H \le 1$. In \cite{G} it was shown that when the support
of a tropical polynomial $f$ is located in $T_R$ then $H(f)\le 1-1/R^n$. \vspace{1mm}

Also one can consider the projective limit $U_{\infty} \subset \RR^{\ZZ^n}$ of $\{U_N\}_{N<\infty}$ with respect to the projections $p_r$ (see above). Then $U_{\infty}$ consists of points $\{u(k_1,\dots,k_n)\, :\, -\infty < k_1,\dots,k_n <\infty\}$
satisfying all the linearizations (\ref{3}). One can view elements of $U_{\infty}$ as $n$-dimensional generalization of tropical recurrent sequences \cite{G}.

For any $N$ and a tropical polynomial $g=\min_{1\le j\le l} \{h_j+\sum_{1\le i\le n} b_{j,i}X_i\}\in rad(I)$ from the radical of $I$ consider its linearization $\min_{1\le j\le l} \{h_j+u(b_{j,1},\dots,b_{j,n})\}$, provided that $0\le b_{j,1},\dots,b_{j,n} < N, \, 1\le j\le l$. Denote by $W_N\subset \RR^{N^n}$ the set of points satisfying these linearizations for all $g\in rad(I)$.

In \cite{GV} an example of intersection of an infinite number of tropical linear prevarieties was produced being not a tropical prevariety. Therefore, it is not clear apriori whether $W_N$ is a tropical linear prevariety (or even a semi-algebraic set). Nevertheless, we define $\dim (W_N)$ as the
minimum of dimensions of tropical linear prevarieties containing $W_N$.

One can prove the existence of the limit
\begin{eqnarray}\label{15}
H(rad(I)):=\lim_{N\to \infty} \dim(W_N)/N^n =\inf \dim (W_N)/N^n
\end{eqnarray}
slightly modifying the above argument which was used to prove the existence of the limit in (\ref{5}). Indeed, consider a partition of
$n$-dimensional grid $T_N$ into subgrids with sides $q_1,\dots,q_R$, respectively. Let $W_{q_r} \subset D_r,\, 1\le r\le R$
where for suitable tropical linear prevarieties $D_r$ hold $\dim(W_{q_r})=\dim (D_r)$. Then $W_N\subset D_{q_1} \times
\cdots \times D_{q_R}$, hence $\dim (W_N)\le \dim(W_{q_1}) + \cdots + \dim (W_{q_R})$ which entails as above (cf. (\ref{5})) the existence of the limit in (\ref{15}).

Similar to $U_{\infty}$ one can consider the projective limit $W_{\infty} \subset \RR^{\ZZ^n}$ of $\{W_N\}_{N<\infty}$.

\begin{remark}
i) Is $W_N$ a tropical linear prevariety? \vspace{2mm}

ii) How to describe the radical of a tropical ideal and $W_N,\, W_{\infty}$ explicitly?
\end{remark}

\section{Vanishing of entropy of  radical of ideal of a tropical curve}\label{two}

In this section we prove that in case of two variables ($n=2$) if all the coefficients of $f$ vanish (i.~e. $a_j=0,\, 1\le j\le m$,
see (\ref{1})) then $H(rad(f))=0$. In  section~\ref{three} an example is exhibited of a non-radical tropical ideal with a positive entropy.

\begin{theorem}\label{curve}
For a tropical bivariate polynomial 
\begin{eqnarray}\label{6}
f=\min_{1\le j\le m} \{t_{j,1}X+t_{j,2}Y\}
\end{eqnarray}
the entropy of its radical $H(rad(f))=0$ equals zero.
\end{theorem} 

{\bf Proof}. Consider Newton polygon of $f$ being the convex hull $P\subset \RR^2$ of the points $(t_{1,1},\, t_{1,2}),\dots, (t_{m,1},\, t_{m,2})$ (see (\ref{6})). Let $P$ be $e$-gon. Then Newton polyhedron $\cal P$ of $f$ is an infinite cylinder in $\RR^3$ (with the coordinates $X,\, Y,\, Z$) such that ${\cal P}= \{(x,\, y,\, z)\, : \, (x,\, y)\in P,\, z\ge 0\}$. The tropical prevariety $V(f)\subset \RR^2$ being moreover, a tropical variety \cite{EKL} (a tropical curve) consists of supporting planes to $\cal P$ in $\RR^3$ (not containing lines parallel to the axis $Z$) which intersect $\cal P$ in at least two points, and thereby, contain an edge of $\cal P$. In fact, $V(f)$ consists of a vertex which corresponds to the face $P$ of $\cal P$ together with $e$ rays emanating from the vertex which correspond to the edges of $\cal P$ (being simultaneously the edges of $P$).

 Observe that if Newton polygon of a tropical polynomial with zero coefficients 
\begin{eqnarray}\label{2}
g=\min_{1\le j\le r} \{b_{j,1}X+b_{j,2}Y\}
\end{eqnarray}
whose Newton polygon contains $e$ edges parallel to the edges of $P$ (and perhaps, in addition, other edges) then $g\in rad(f)$.
This differs the tropical situation from  the classical commutative algebra with respect to Hilbert's strong Nullstellensatz (providing the structure of the radical of an ideal). \vspace{1mm}

To prove the Theorem it suffices to verify that $\dim (W_N)=o(N^2)$. Moreover, we'll prove that $\dim (W_N)=O(N)$.

Denote by $W_N^{(d)} \subset W_N$ the subset of points $\{w(k_1,k_2)\, :\, 0\le k_1,k_2< N\} \subset W_N$ such that among the values of $w(k_1,k_2),\,  0\le k_1,k_2< N$ there are at most $d$ different. We'll prove that $W_N^{(d)}=W_N^{(d-1)}$ for $d>c_0\cdot N$ for an appropriate

\noindent constant $c_0>0$, in other words, $W_N=W_N^{(d)}$ for $d=\lfloor c_0\cdot N \rfloor$. Then $W_N^{(d)}$ is contained in a tropical linear prevariety (i.~e. in a finite union of polyhedra) 
$$\{w(k_1,k_2)\, :\, 0\le k_1,k_2< N,\, \mbox{among the values of}\, w(k_1,k_2)\, \mbox{are at most}\, d\, \mbox{different}\}$$

\noindent Since the dimension of the latter set does not exceed $d$, one can talk about upper bounds on $\dim (W_N)$ (see section~\ref{one}), and thereby we'll prove an upper bound $\dim (W_N) =O(N)$ (which suffices for the proof of the Theorem).
\vspace{1mm}

 Fix a point $w_0=\{w(k_1,k_2)\, :\, 0\le k_1,k_2< N\} \in W_N^{(d_0)}$ for some $d_0$. We describe a recursive process in the course of which it modifies (more precisely, shrinks) a polygon $Q$. As a base of recursion we take as $Q$ the square $\{(x,y)\, :\, 0\le x,y< N\}$. At every step of the recursion take an edge $E$ of (Newton) polygon $P$. Denote by $L:=L(E)$ the line containing edge $E$. First we move $L$ parallel to itself outwards $P$ until we reach a line $L_0:=L_0(E)$ such that $P$ and the (current) $Q$ lie on the same side of $L_0$. Then we move $L_0$ parallel to itself in the direction towards $P$ (let us call it for definiteness, the inwards direction) until we reach a line $L_1:=L_1(E)$ which for the first time contains an integer point from $Q$. Therefore, all the integer points from $Q$ are situated in the half-plane bounded by $L_1$. If $w_0(k_1,k_2)$ takes at all the integer points $(k_1,k_2)\in L_1\cap Q$ at most two different values, then we move further $L_1$  parallel to itself in the inwards direction until the resulting line $L_2$ reaches an integer point of $Q$ (unless there are no other integer points in $Q$, and the process terminates). Thus, all the integer points from $Q$ lie either on $L_1$ or in the half-plane $S$ bounded by $L_2$. For the next step of the recursion we shrink $Q$ intersecting it with $S$.

Now alternatively, we suppose that for each edge $E$ of $P$ the constructed above line $L_1(E) \cap Q$ intersected with $Q$ contains at least three integer points having different values of $w_0$. Choose a triple of such points for each edge (the triples may intersect). Denote by $A$ the union over all the edges of $P$ of these triples of points. Among the points from $A$ choose a point $(k_1,k_2)$ with the minimal value of $w_0(k_1,k_2)$. Then for each edge $E$ of $P$ choose two integer points from $L_1(E)\cap A$ having values of $w_0$ greater than $w_0(k_1,k_2)$. The set of chosen points (including point $(k_1,k_2)$) denote by $B\subset A$.

Consider a tropical polynomial $g:=\min_{(b_1,b_2)\in B} \{b_1X+b_2Y\}$. Then $g\in rad(f)$ because for each edge $E$ of   $P$ Newton polygon of $g$ contains an edge parallel to $E$. We claim that point $w_0\in W_N^{(d_0)} \subset W_N$ does not satisfy the linearization $\min_{(b_1,b_2)\in B} \{w(b_1,b_2)\}$ of $g$ (cf. (\ref{3})). Assume the contrary. We have $w_0(k_1,k_2)= \min_{(b_1,b_2)\in B} \{w_0(b_1,b_2)\}$. There is a point $B\ni (a_1,a_2) \neq (k_1,k_2)$ such that $w_0(k_1,k_2)=w_0(a_1,a_2)$. This contradicts to the choice of $B$.

Thus, the described above recursive process terminates only with exhausting all the integer points of grid $T_N$. Observe that there are at most $O(N)$ steps of the process because consecutive parallel lines $L_1,\, L_2$ (see above) contain integer points, and taking into account that there are at most $O(N)$ lines parallel to $E$ and intersecting $T_N$ (a constant hidden in big "O" is determined by the denominators of the slopes of the edges of $P$, so depends on $t_{j,1},\, t_{j,2},\, 1\le j\le m$). At each step of the process at most two different values of $w_0$ occur. Hence $w_0$ has at most $O(N)$ different values, i.~e. $d_0=O(N)$  which completes the proof of the Theorem~\ref{curve}. $\Box$ 

\begin{remark}
We have shown in the proof that $\dim (W_N)=O(N)$. Does there exist the limit $\lim_{N\to \infty} \dim (W_N)/N$? If it were the case this limit would play a role of a tropical version of the leading coefficient of Hilbert's polynomial of the radical ideal $rad(f)$.  
\end{remark}

\section{Bounds on entropy of tropical polynomial $\min\{0,\, X,\, Y,\, X+Y\}$}\label{three}

In this section we provide more precise bounds on dimensions $\dim (W_N)$ for tropical polynomial $f:= \min\{0,\, X,\, Y,\, X+Y\}$ and show that $H(f)>0$ (unlike $H(rad(f))=0$ according to the Theorem~\ref{curve}). Note that Newton polygon $P$ of $f$ is $1\times 1$ square $\{(x,y)\, :\, 0\le x,y\le 1\}$.

Following the recursive process described in the proof of the Theorem~\ref{curve} one can observe that in case of $f$ all the intermediate polygons in the course of the process are rectangles with horizontal and vertical sides. Moreover, horizontal sides of rectangles can be moved at most $N$ times (as well as vertical sides). Thus, the recursive process runs in at most $2N$ steps. At every step at most two values of $w_0$ occur. Hence $\dim (W_N)\le 4N$. \vspace{1mm}

On the other hand, to establish a lower bound on $\dim (W_N)$, consider the following point $\{w_{\infty}(x,y)=c(x)\, :\, x,y \in \ZZ\} \in W_{\infty} \subset \RR^{\ZZ^2}$ from the infinite-dimensional space on $\ZZ^2$ for an arbitrary concave function $c$. One can verify that $w_{\infty}$ satisfies all the linearizations of $rad(f)$. Thus, restricting $w_{\infty}$ on $N\times N$ grid $T_N\subset \ZZ^2$ we conclude that $W_N$ contains a tropical linear prevariety of dimension $N$, therefore $\dim (W_N) \ge N$.

Another example is $\{w_{\infty}(x,y)=0,\, w_{\infty}(x,x)\ge 0\, :\, x,y\in \ZZ,\, x\neq y\}  \subset W_{\infty} \subset \RR^{\ZZ^2}$, providing a lower bound $\dim (W_N) \ge N +1$. \vspace{1mm}

Finally, we show that $H(f)>0$, in other words, the condition of a tropical ideal to be radical in the Theorem~\ref{curve}, is essential. Consider the following set
$$\{u_{\infty}(2x,y)=0,\, u_{\infty}(2x+1,y)\ge 0\, :\, x,y\in \ZZ \} \subset U_{\infty} \subset \RR^{\ZZ^2}.$$
\noindent Restricting $u_{\infty}$ on grid $T_N$  (obtaining a subset of $U_N \subset \RR^{N^2}$) we get that $\dim (U_N)\ge N\cdot \lfloor N/2 \rfloor$, therefore $H(f)\ge 1/2$.

More exotic examples of  points from $U_{\infty}$ one can find in \cite{GP}. 

\section{Vanishing entropy of radical ideal of zero-dimensional tropical prevariety}\label{four}

Let $V\subset \RR^n$ be a tropical prevariety consisting of a finite number $k$ of points. One can treat elements of $V$ as hyperplanes in $\RR^{n+1}$ given by linear equations  $$Z=\sum_{1\le j\le n} l_{i,j}X_j:=L_i(X_1,\dots,X_n),\, 1\le i\le k,$$ \noindent respectively.  In this section we prove the following result.

\begin{theorem}\label{zero}
For a zero-dimensional tropical prevariety $V$ the entropy of its radical $H(rad(V))=0$.
\end{theorem}

{\bf Proof}. \subsection{Radical of zero-dimensional tropical prevariety}\label{five}

First we describe a construction of all the elements of the radical ideal $rad(V)$. Let $\RR^{n+1}$ be equipped with the coordinates $X_1,\dots,X_n,Z$, the coordinate $Z$ we call vertical. Consider a polyhedron $Q\subset \RR^{n+1}$ with the facets parallel to the hyperplanes $\{Z=L_1\},\, \dots,\{Z=L_k\}$, respectively, and with any its point $(a_1,\dots,a_n,b_0)\in Q$ containing the vertical ray $\{(a_1,\dots,a_n,b)\, :\, b\ge b_0\} \subset Q$ emanating from $(a_1,\dots,a_n,b_0)$. In addition, we assume that $Q$ has $2n$ vertical facets contained in the hyperplanes $\{X_j=-1\},\, \{X_j=N\},\, 1\le j\le n$, respectively, for some
sufficiently big (varying) $N$. Thus, $Q$ has either compact facets parallel to  $\{Z=L_1\},\, \dots,\{Z=L_k\}$ or vertical ones. 
The projection on the hyperplane $\{Z=0\}$ coincides with the $n$-dimensional cube $\{-1\le X_1,\dots,X_n\le N\}$.

For an integer point $(a_1,\dots,a_n)\in \ZZ^n$ we call a point   $(a_1,\dots,a_n,b)$ for $b\in \RR$ also integer (abusing the language without leading to misunderstanding). Pick on every compact facet of $Q$ at least two integer points with $(a_1,\dots,a_n)\in \{0,\dots,N-1\}^n$, provided that it is possible. In addition, pick a finite set of integer points strictly inside $Q$. Denote by $A\subset \RR^{n+1}$ the set of  all picked points. Then the tropical polynomial 
\begin{eqnarray}\label{10}
\min_{(a_1,\dots,a_n,b)\in A} \{a_1X_1+\cdots+a_nX_n+b\}
\end{eqnarray}
belongs to the radical $rad(V)$. In fact, conversely, any element of $rad(V)$ can be obtained in this way (for an appropriate $N$).

Making a suitable homothety of $Q$ (thereby, changing $N$), one can suppose that every compact facet of $Q$ contains at least two integer points in the inner (so, relatively open) part  of the facet (in other words, not on the boundary of the facet). 

Since apriori it is unclear whether $W_N$ (see section~\ref{one}) is a tropical linear prevariety, or more generally, a semi-algebraic set, our purpose is to produce semi-algebraic sets (moreover, tropical linear prevarieties) $W_N\subset {\cal W}_N \subset \RR^{N^n}$ such that $\dim ( {\cal W}_N)=O(N^{n-1})$.

Fix a point $w=\{(a_1,\dots,a_n,\, w(a_1,\dots,a_n))\, :\, (a_1,\dots,a_n)\in \{0,\dots,N-1\}^n\}\in W_N$ for the time being. Shifting $Q$ in the vertical direction (so, along the axis $Z$) one can assume that no integer point among $\{(a_1,\dots,a_n,\, -w(a_1,\dots,a_n))\, :\, (a_1,\dots,a_n)\in \{0,\dots,N-1\}^n\}$ lies strictly inside $Q$, and on the other hand, at least one of these points lies on the boundary of $Q$.

\subsection{Regular and singular facets of a polyhedron from the radical. Extending a polyhedron}\label{nine} 

Now we'll  describe the process of moving (some of) the facets of $Q$ parallel to themselves outwards $Q$, thereby, extending $Q$. At the beginning of the process all the facets of $Q$ are declared {\it singular}. In the way of the process some of the facets become {\it regular}. In the latter case they remain regular in the whole way of the process, and we don't move them anymore.

Note that $Q$ is defined uniquely by the hyperplanes $\{Z=L_i+c_i\},\, 1\le i\le k$ with $c_1,\dots,c_k \in \RR$ which contain the facets of $Q$. Abusing the language (without leading to misunderstanding) we denote the facets by their corresponding hyperplanes $\{Z=L_i+c_i\},\, 1\le i\le k$.

Thus, at each step of the process we fix  an arbitrary current singular facet $\{Z=L_i+c_i\}$ for some $1\le i\le k$, keep the current facets $\{Z=L_s+c_s\},\, 1\le s\le k, s\neq i$ and replace the facet  $\{Z=L_i+c_i\}$ by $\{Z=L_i+c'_i\}$ for an appropriate $c'_i \le c_i$. The resulting extended polyhedron we denote by $Q_i$. Still, $Q_i$ has vertical unbounded facets contained in the hyperplanes $\{X_j=-1\},\, \{X_j=N\},\, 1\le j\le n$, respectively (cf. section~\ref{five}).

Below we consider polyhedra $Q_{i,c}$ with the fixed facets $\{Z=L_s+c_s\},\, 1\le s\le k, s\neq i$ and $\{Z=L_i+c\}$ with varying $c$. Take the maximal $c$ such that there exists (provided that it does exist) an integer point $(a_1,\dots,a_n) \in \{0,\dots,N-1\}^n$ such that $-w(a_1,\dots,a_n)=L_i(a_1,\dots,a_n)+c$ and in addition, the point $(a_1,\dots,a_n,-w(a_1,\dots,a_n))$ belongs to the inner part of the facet $\{Z=L_i+c\}$ of the polyhedron $Q_{i,c}$, i.~e. $-w(a_1,\dots,a_n) > L_s(a_1,\dots,a_n)+c_s$ for $1\le s\le k, s\neq i$. Assume that there exists another integer point $(a'_1,\dots,a'_n) \in \{0,\dots,N-1\}^n,\, (a'_1,\dots,a'_n) \neq (a_1,\dots,a_n)$ for which the point $(a'_1,\dots,a'_n,L_i(a'_1,\dots,a'_n)+c)$ also belongs to the inner part of the facet $\{Z=L_i+c\}$ of $Q_{i,c}$. Then we put $c'_i:=c,\, Q_i:=Q_{i,c}$ and declare the moved facet $\{Z=L_i+c'_i\}$ to be {\it regular} in $Q_i$.

Otherwise, if the conditions of the previous paragraph are not fulfilled, we take the minimal $c'$ such that the facet $\{Z=L_i+c'\}$ of $Q_{i,c'}$ contains at least two integer points of the form $(a_1,\dots,a_n,L_i(a_1,\dots,a_n)+c')$ with $(a_1,\dots,a_n) \in \{0,\dots,N-1\}^n$, i.~e. 
$$L_i(a_1,\dots,a_n)+c'\ge L_s(a_1,\dots,a_n)+c_s,\, 1\le s\le k, s\neq i.$$ 
\noindent Note that $c'\ge c$ (provided that $c$ from the previous paragraph does exist, while no point $(a'_1,\dots,a'_n)$ satisfying the properties from the previous paragraph does exist). In the case under consideration we put $c'_i:=c',\, Q_i:= Q_{i,c'}$, the (moved) facet $\{Z=L_i+c'\}$ is declared to be {\it singular}.

Thus, we have yielded an (extended) polyhedron $Q_i$. Observe that still no integer point among $(a_1,\dots,a_n,-w(a_1,\dots,a_n))$ lies strictly inside $Q_i$. Also observe that if an integer point $(a_1,\dots,a_n,b)$ lies in (respectively, the inner part of) a facet $\{Z=L_s+c_s\},\, 1\le s\le k, s\neq i$  of $Q$ then $(a_1,\dots,a_n,b)$ lies in (respectively, the inner part of) the facet $\{Z=L_s+c_s\}$ of $Q_i$.

We apply the described process consecutively to singular facets (in an arbitrary order), while it is possible, extending a current polyhedron $Q$, the resulting polyhedron we denote by $Q(\infty)$. If the process is infinite we show that still it tends to a limit polyhedron $Q(\infty)$. For each $1\le i\le k$  if a current facet $\{Z=L_i+c_i(t)\}$ at the moment $t$ is moved infinite number of times  then there exists an integer point $(a_1,\dots,a_n) \in \{0,\dots,N-1\}^n$ such that for an infinite number of moments $t$ a point $(a_1,\dots,a_n,L_i(a_1,\dots,a_n)+c_i(t))$ lies on the facet $\{Z=L_i+c_i(t)\}$ (moreover, there exist at least two such points). Therefore, $L_i(a_1,\dots,a_n)+c_i(t)\ge -w(a_1,\dots,a_n)$ because the point $(a_1,\dots,a_n,-w(a_1,\dots,a_n))$ lies outside of the inner part of the current polyhedron. Hence there exists the limit $c_i(\infty)=\lim _{t\to \infty} c_i(t)$ (recall that $c_i(t)$ decrease with growing $t$). Observe that the point $(a_1,\dots,a_n,L_i(a_1,\dots,a_n)+c_i(\infty))$ belongs to the facet $\{Z=L_i+c_i(\infty)\}$ of the limit polyhedron $Q(\infty)$ (thus, there exist at least two such points). Indeed, 
$$L_i(a_1,\dots,a_n)+c_i(t)\ge \max _{1\le s\le k,\, s\neq i} \{L_s(a_1,\dots,a_n)+c_s(t)\}$$ 
\noindent for an infinite number of moments $t$, and this inequality holds also after taking the limit.

Let a facet $\{Z=L_i+c_i(\infty)\}$ of $Q(\infty)$ be singular for some $1\le i\le k$ (it means that the facet $\{Z=L_i+c_i(t)\}$ was  singular for any $t$). We prove the following bound.

\begin{lemma}\label{bound}
The number of integer points $(a_1,\dots,a_n) \in \{0,\dots,N-1\}^n$ such that the point $(a_1,\dots,a_n,L_i(a_1,\dots,a_n)+c_i(\infty))$ belongs to a singular facet $\{Z=L_i+c_i(\infty)\}$ of $Q(\infty)$, does not exceed $O(N^{n-1})$.
\end{lemma}

{\bf Proof of the Lemma}. To prove this bound observe that all the integer points (except of maybe one point) of the form $(a_1,\dots,a_n,L_i(a_1,\dots,a_n)+c_i(\infty)),\, (a_1,\dots,a_n) \in \{0,\dots,N-1\}^n$ which belong to the facet $\{Z=L_i+c_i(\infty)\}$ are located on the boundary of this facet. Indeed, otherwise, if there were two points not in the boundary, one could further move the facet $\{Z=L_i+c_i(\infty)\}$, i.~e. decrease $c_i(\infty)$. Thus, all the points under consideration (except of maybe one point) lie in the intersection of a pair of hyperplanes $\{Z=L_i+c_i(\infty)\}$ and $\{Z=L_s+c_s(\infty)\}$ for a certain $1\le s\le k,\, s\neq i$. Every such intersection is a plane of dimension $n-1$, and its projection on the hyperplane $\{Z=0\}$ contains at most of $O(N^{n-1})$ points among $\{0,\dots,N-1\}^n$. This proves the required bound $O(N^{n-1})$ on the number of integer points in the facet. $\Box$ 

\subsection{Slight shrinking of a polyhedron from the radical}\label{six}

Now we choose a sufficiently small $\varepsilon >0$ and replace each singular facet $\{Z=L_i+c_i(\infty)\}$ of $Q(\infty)$ by the facet $\{Z=L_i+c_i(\infty)+\varepsilon \}$ (and keeping the regular facets, see section~\ref{nine}). Thereby, we shrink the polyhedron  $Q(\infty)$, the resulting polyhedron denote by $Q^{(\varepsilon)}(\infty)$. We require (choosing $\varepsilon$ sufficiently small) that any integer point of the form  $(a_1,\dots,a_n,L_r(a_1,\dots,a_n)+c_r(\infty))$ lying in the inner part of a regular facet $\{Z=L_r+c_r(\infty)\}$ of $Q(\infty)$ still lies in the inner part  of the facet $\{Z=L_r+c_r(\infty)\}$ of $Q^{(\varepsilon)}(\infty)$ (recall, see section~\ref{nine}, that there exist at least two such points).

Note that for any integer point of the form $(a_1,\dots,a_n,L_i(a_1,\dots,a_n)+c_i(\infty)),\, (a_1,\dots,a_n) \in \{0,\dots,N-1\}^n$ which lies in a singular facet $\{Z=L_i+c_i(\infty)\}$ of the polyhedron $Q(\infty)$ the point $(a_1,\dots,a_n,L_i(a_1,\dots,a_n)+c_i(\infty)+\varepsilon)$ lies in the facet $\{Z=L_i+c_i(\infty)+\varepsilon\}$ of the polyhedron $Q^{(\varepsilon)}(\infty)$, therefore this facet $\{Z=L_i+c_i(\infty)+\varepsilon\}$ contains at least two such integer points, see section~\ref{nine} (and at most $O(N^{n-1})$ integer points due to Lemma~\ref{bound}). 

Assume that there is a singular facet $\{Z=L_r+c_r(\infty)\}$ of $Q(\infty)$ such that each integer point (except of maybe one point) of the form  $(a_1,\dots,a_n,L_r(a_1,\dots,a_n)+c_r(\infty)),\, (a_1,\dots,a_n) \in \{0,\dots,N-1\}^n$ lying in the facet $\{Z=L_r+c_r(\infty)\}$ of the polyhedron $Q^{(\varepsilon)}(\infty)$ coincides with the point $(a_1,\dots,a_n,-w(a_1,\dots,a_n))$. In other words, the latter point $(a_1,\dots,a_n,-w(a_1,\dots,a_n))$ lies in the facet $\{Z=L_r+c_r(\infty)\}$ of the polyhedron $Q^{(\varepsilon)}(\infty)$, as well as this point lies in the facet $\{Z=L_r+c_r(\infty)\}$ of the polyhedron $Q(\infty)$. Then we further shrink the polyhedron $Q^{(\varepsilon)}(\infty)$ replacing the facet $\{Z=L_r+c_r(\infty)\}$ by $\{Z=L_r+c_r(\infty)+\varepsilon_0\}$ (and keeping all other facets of $Q^{(\varepsilon)}(\infty)$), the resulting polyhedron we denote by $Q_r^{(\varepsilon)}(\infty)$. We choose $\varepsilon_0 >0$ sufficiently small such that any integer point of the form $(a_1,\dots,a_n,b),\, (a_1,\dots,a_n) \in \{0,\dots,N-1\}^n$ lying in the inner part of a (either regular or singular) facet (different from $\{Z=L_r+c_r(\infty)\}$) of the polyhedron $Q^{(\varepsilon)}(\infty)$ still lies in the inner part of the corresponding facet of the resulting polyhedron $Q_r^{(\varepsilon)}(\infty)$.

Moreover, we require (choosing $\varepsilon_0 >0$ sufficiently small) that any point of the form $(a_1,\dots,a_n,L_i(a_1,\dots,a_n)+c_i(\infty)+\varepsilon)$ lying in the facet $\{Z=L_i+c_i(\infty)+\varepsilon\}$ of the polyhedron $Q^{(\varepsilon)}(\infty)$ for a singular facet $\{Z=L_i+c_i(\infty)\}$ of the polyhedron $Q(\infty)$ still lies in the facet $\{Z=L_i+c_i(\infty)+\varepsilon\}$ of the polyhedron $Q_r^{(\varepsilon)}(\infty)$ (cf. section~\ref{nine}). This choice is possible because there are no common integer points of the form $(a_1,\dots,a_n,b),\, (a_1,\dots,a_n) \in \{0,\dots,N-1\}^n$ of the facets $\{Z=L_r+c_r(\infty)\}$ and $\{Z=L_i+c_i(\infty)+\varepsilon\}$ of the polyhedron $Q^{(\varepsilon)}(\infty)$ due to the choice of $\varepsilon$.
Therefore, each facet of the form $\{Z=L_i+c_i(\infty)+\varepsilon\}$ of the polyhedron $Q_r^{(\varepsilon)}(\infty)$ still contains at least two integer points (and at most $O(N^{n-1})$ ones due to Lemma~\ref{bound}).

If the polyhedron $Q_r^{(\varepsilon)}(\infty)$ contains another regular facet of the polyhedron $Q(\infty)$ satisfying the same properties as $\{Z=L_r+c_r(\infty)\}$ from the previous two paragraphs, we shrink further the polyhedron $Q_r^{(\varepsilon)}(\infty)$, and continuer this process, while it is possible. Note that during this process each regular facet of $Q(\infty)$ can be taken at most once, so the process terminates after at most of $k$ steps. At the end we obtain a polyhedron which we denote by $Q_{fin}^{(\varepsilon)}(\infty)$. \vspace{1mm}

First assume that the polyhedron $Q(\infty)$ contains a regular facet $\{Z=L_r+c_r(\infty)\}$ which was not moved in the process of constructing $Q_{fin}^{(\varepsilon)}(\infty)$ from $Q^{(\varepsilon)}(\infty)$. Then the inner part of the facet $\{Z=L_r+c_r(\infty)\}$ of the polyhedron $Q_{fin}^{(\varepsilon)}(\infty)$ contains at least two integer points, among which there is a point 
\begin{eqnarray}\label{11}
(a_1^{(0)},\dots,a_n^{(0)},(L_r(a_1^{(0)},\dots,a_n^{(0)})+c_r(\infty)=-w(a_1^{(0)},\dots,a_n^{(0)}))) 
\end{eqnarray}
for some $(a_1^{(0)},\dots,a_n^{(0)}) \in \{0,\dots,N-1\}^n$ (see section~\ref{nine}). On the other hand, the facet $\{Z=L_r+c_r(\infty)\}$ contains at least two integer points 
$$(a_1^{'},\dots,a_n^{'},(L_r(a_1^{'},\dots,a_n^{'}+c_r(\infty))),\, (a_1^{''},\dots,a_n^{''},(L_r(a_1^{''},\dots,a_n^{''}+c_r(\infty)))$$
\noindent for some points $(a_1^{'},\dots,a_n^{'}),\, (a_1^{''},\dots,a_n^{''}) \in \{0,\dots,N-1\}^n$ such that
$$-w(a_1^{'},\dots,a_n^{'})<L_r(a_1^{'},\dots,a_n^{'})+c_r(\infty),\, -w(a_1^{''},\dots,a_n^{''})<
L_r(a_1^{''},\dots,a_n^{''})+c_r(\infty),$$
\noindent i.~e. the points $(a_1^{'},\dots,a_n^{'},-w(a_1^{'},\dots,a_n^{'})),\, (a_1^{''},\dots,a_n^{''},-w(a_1^{''},\dots,a_n^{''}))$ are located outside the polyhedron $Q_{fin}^{(\varepsilon)}(\infty)$.

Now we produce a family of points $A\subset \RR^{n+1}$. First, we put two points 
\begin{eqnarray}\label{12}
(a_1^{(0)},\dots,a_n^{(0)},-w(a_1^{(0)},\dots,a_n^{(0)})),\, (a_1^{'},\dots,a_n^{'},L_r(a_1^{'},\dots,a_n^{'})+c_r(\infty))
\end{eqnarray}
 in $A$. 

Then for each facet $\{Z=L_s+c_s(\infty)\}$ of the polyhedron $Q_{fin}^{(\varepsilon)}(\infty)$ which was not moved in the process of constructing $Q_{fin}^{(\varepsilon)}(\infty)$ from $Q^{(\varepsilon)}(\infty)$, and the facet $\{Z=L_s+c_s(\infty)\}$ is a regular one of the polyhedron $Q(\infty)$, take two integer points 
$$a^{(1)}:=(a_1^{(1)},\dots,a_n^{(1)},L_s(a_1^{(1)},\dots,a_n^{(1)})+c_s(\infty)),$$ $$a^{(2)}:=(a_1^{(2)},\dots,a_n^{(2)},L_s(a_1^{(2)},\dots,a_n^{(2)})+c_s(\infty))$$
\noindent lying in the facet $\{Z=L_s+c_s(\infty)\}$ of the polyhedron $Q_{fin}^{(\varepsilon)}(\infty)$ such that 
$$L_s(a_1^{(1)},\dots,a_n^{(1)})+c_s(\infty)>-w(a_1^{(1)},\dots,a_n^{(1)}),$$
$$L_s(a_1^{(2)},\dots,a_n^{(2)})+c_s(\infty)>-w(a_1^{(2)},\dots,a_n^{(2)}).$$
\noindent Such points exist because the facet $\{Z=L_s+c_s(\infty)\}$ was not moved. We add the points $a^{(1)},\, a^{(2)}$ to $A$. 

Finally, for every other facet of  the polyhedron $Q_{fin}^{(\varepsilon)}(\infty)$ we take two arbitrary integer points lying in this facet and add them to $A$ (note that some of the added points can be common for different facets). Observe that for each of these points $(a_1,\dots,a_n,b)$ it holds $b>-w(a_1,\dots,a_n)$ (due to shrinking $Q(\infty)$ while constructing $Q^{(\varepsilon)}(\infty)$ and subsequently further shrinking while constructing $Q_{fin}^{(\varepsilon)}(\infty)$).

Consider a tropical polynomial 
$$f:=\min_{(a_1,\dots,a_n,b)\in A} \{a_1X_1+\cdots +a_nX_n+b\}.$$
\noindent Then $f\in rad(V)$ (cf. (\ref{10})). On the other hand, the point $-w$ is not a tropical solution of the linearization 
\begin{eqnarray}\label{13}
\min_{(a_1,\dots,a_n,b)\in A} \{w(a_1,\dots,a_n)+b\}
\end{eqnarray}
 of $f$ (cf. (\ref{3})). Indeed, at the point $(a_1^{(0)},\dots,a_n^{(0)},-w(a_1^{(0)},\dots,a_n^{(0)})) \in A$ (see (\ref{11}), (\ref{12})) the value of the linearization (\ref{13}) equals $0$. At any other point  $(a_1,\dots,a_n,b)\in A$ the value of the linearization is positive due to the construction of $A$. Whence we get a contradiction with the assumption that the polyhedron $Q(\infty)$ contains a (regular) facet which was not moved in the process of constructing the polyhedron $Q_{fin}^{(\varepsilon)}(\infty)$.

\subsection{Tropical solution of the linearization of the radical almost lies on the boundary of a polyhedron}\label{eight}

Now we bound from below the number of integer points $(a_1,\dots,a_n) \in \{0,\dots,N-1\}^n$ such that the point $(a_1,\dots,a_n,-w(a_1,\dots,a_n))$ lies in a boundary of the polyhedron $Q(\infty)$. For each integer point $(a_1,\dots,a_n) \in \{0,\dots,N-1\}^n$ there is a unique point $(a_1,\dots,a_n,b)$ lying in the boundary of $Q(\infty)$. 

If $(a_1,\dots,a_n,b)$ belongs to a regular facet  $\{Z=L_r+c_r(\infty)\}$ of $Q(\infty)$ (cf. section~\ref{nine}) and does not belong to any singular facet of $Q(\infty)$ then $b=-w(a_1,\dots,a_n)$ holds for all such points $(a_1,\dots,a_n)$ (except, perhaps, of one point) for a given regular facet $\{Z=L_r+c_r(\infty)\}$ (due to the construction of the polyhedron $Q_{fin}^{(\varepsilon)}(\infty)$). 

Otherwise, if the point $(a_1,\dots,a_n,b)$ belongs to a singular facet of $Q(\infty)$ then according to Lemma~\ref{bound} there are at most $O(N^{n-1})$ integer points $(a_1,\dots,a_n) \in \{0,\dots,N-1\}^n$ such that the point $(a_1,\dots,a_n,b)$ belongs to a singular facet of $Q(\infty)$.

Thus, totally for at most of $d_0N^{n-1}$ (for an appropriate constant $d_0>0$) integer points $(a_1,\dots,a_n)$ among $\{0,\dots,N-1\}^n$ the point $(a_1,\dots,a_n,-w(a_1,\dots,a_n))$ does not belong to the boundary of the polyhedron $Q(\infty)$ (taking into account that $Q(\infty)$ has a constant number $k$ of compact facets).

Therefore the point 
$$w^{(-)}:=\{(a_1,\dots,a_n,-w(a_1,\dots,a_n))\, :\, (a_1,\dots,a_n) \in \{0,\dots,N-1\}^n\} \in \RR^{N^n}$$
\noindent belongs to the defined below semi-algebraic set (in fact, it is a finite union of polyhedra, so a tropical linear prevariety) ${\cal W}_N^{(-)} \subset \RR^{N^n}$ (for the sake of convenience of notations we replace the point $w$ by $w^{(-)}$). Take an arbitrary subset $D$ of $\lfloor d_0N^{n-1} \rfloor$ points in $ \{0,\dots,N-1\}^n$. Then the semi-algebraic set ${\cal W}_N^{(D)}$ consists of all the points $\{(a_1,\dots,a_n,b) \, :\, (a_1,\dots,a_n) \in \{0,\dots,N-1\}^n\} \in \RR^{N^n}$ for which there exists a vector $(c_1,\dots,c_k) \in \RR^k$ such that each point $(a_1,\dots,a_n,b)$ when $(a_1,\dots,a_n) \not \in D$, lies in the boundary of the convex polyhedron $\{Z\ge L_i+c_i,\, 1\le i\le k\}$. We define ${\cal W}_N^{(-)}$ as the union of sets ${\cal W}_N^{(D)}$ over all the subsets $|D|=\lfloor d_0N^{n-1} \rfloor$. Clearly, $\dim({\cal W}_N^{(D)})\le d_0N^{n-1}+k$. Thus, $\dim({\cal W}_N^{(-)}) =O(N^{n-1})$ which completes the proof of the Theorem~\ref{zero}. $\Box$

 \vspace{2mm}

{\bf Acknowledgements}. The author is grateful to the grant RSF 16-11-10075 and
 to MCCME for inspiring atmosphere.

\end{document}